\documentclass[12pt,anrf]{article}

\newtheorem{thm}{Theorem}[section]

\newtheorem{lem}[thm]{Lemma}

\begin{document}

\title{{\bf A New Radial Function}}
\author{{\bf Lin-Tian Luh\thanks{This work was supported by NSC93-2115-M126-004.}}}

\maketitle

\bigskip
{\bf Abstract.}  In the field of radial basis functions mathematicians have been endeavouring to find infinitely differentiable and compactly supported radial functions. This kind of functions is extremely important. One of the reasons is that its error bound will converge very fast. However there is hitherto no such function which can be expressed in a simple form. This is a famous question. The purpose of this paper is to answer this question.\\
\\
{\bf keywords}: radial basis function, compactly supported function, positive definite function.

\pagenumbering{arabic}
\setcounter{section}{0}

\section{Introduction}

In this paper we raise a new function
\begin{equation}
  \Phi(x):=\left\{ \begin{array}{ll}
                     e^{-\alpha(1+tan\frac{\pi}{2}\| x\| ^{2})^{2}} & \mbox{if $\| x\| <1$}\\
                     0                                              & \mbox{if $\| x\| \geq 1$}
                   \end{array}
           \right.                                            
\end{equation}
where $\alpha > 0,\ x\in R^{n}$ and $\| x\| $ denotes the Euclidean norm of $x$. The infinite differentiability of $\Phi$ can be shown by routine investigation, and we omit it. What's difficult is to show its positive definiteness. In order to make it useful in scattered data interpolation, this problem has to be overcome. This is the main theme of the next section.

\section{Positive Definiteness}

It's well known that positive definiteness is just conditional positive definiteness of order zero. Here we adopt the definition of \cite{Dy} for conditional positive definiteness. Before showing that the map defined in (1) is positive definite, let's define
\begin{equation}
  \phi(t):=\left\{ \begin{array}{ll}
                     e^{-\alpha(1+tan\frac{\pi}{2}t^{2})^{2}} & \mbox{if $0\leq t<1$}\\
                     0                                        & \mbox{if $1\leq t.$}
                   \end{array}
           \right.   
\end{equation}
In other words, $\phi(\|x\|)=\Phi(x)$ for $x\in R^{n}$, where $\alpha>0$ is a constant. In order to make this paper easier to understand, let's analyze the derivatives of $\phi(\sqrt{t})$ of order up to four.

It's easily seen that with $u=tan\frac{\pi}{2}t$,
\begin{itemize}
  \item $\frac{d}{dt}\phi(\sqrt{t})=e^{-\alpha(1+u)^{2}}(-\alpha)\frac{\pi}{2}(1+u^{2})[2u+2]$
  \item $\frac{d^{2}}{dt^{2}}\phi(\sqrt{t})=e^{-\alpha(1+u)^{2}}\alpha(\frac{\pi}{2})^{2}(1+u^{2})[4\alpha(1+u)^{2}(1+u^{2})-2(1+u^{2})-4(1+u)u]$
  \item $\frac{d^{3}}{dt^{3}}\phi(\sqrt{t})=e^{-\alpha(1+u)^{2}}(-\alpha)(\frac{\pi}{2})^{3}(1+u^{2})[8\alpha^{2}(1+u)^{3}(1+u^{2})^{2}-12\alpha(1+u)(1+u^{2})^{2}-24\alpha u(1+u)^{2}(1+u^{2})+8(1+u^{2})u+4(1+2u)(1+u^{2})+8(1+u)u^{2}]$
  \item $\frac{d^{4}}{dt^{4}}\phi(\sqrt{t})=e^{-\alpha(1+u)^{2}}\alpha(\frac{\pi}{2})^{4}(1+u^{2})[16\alpha^{3}(1+u)^{4}(1+u^{2})^{3}-48\alpha^{2}(1+u)^{2}(1+u^{2})^{3}-96\alpha^{2}(1+u)^{3}u(1+u^{2})^{2}+136\alpha(1+u)u(1+u^{2})^{2}+8\alpha(1+u)(1+2u)(1+u^{2})^{2}+64\alpha(1+u)^{2}u^{2}(1+u^{2})+12\alpha(1+u^{2})^{3}+24\alpha(1+u)^{2}(1+3u^{2})(1+u^{2})-104u^{2}(1+u^{2})-32u(1+u^{2})-16(1+u^{2})-16u^{3}(1+u)].$
\end{itemize}
By induction, we find
\begin{equation}
  \frac{d^{j}}{dt^{j}}\phi(\sqrt{t})=e^{-\alpha(1+u)^{2}}(-1)^{j}\alpha(\frac{\pi}{2})^{j}(1+u^{2})F_{j}(\alpha, u).
\end{equation}
where $F_{j}(\alpha, u)$ is a polynomial of $\alpha,\ u$ with integer coefficients. Note that $$F_{j}(\alpha,u)=a_{j}\alpha^{j-1}u^{3(j-1)+1}+\sum_{k=0}^{j-2}\sum_{l=0}^{3j-3}a_{kl}\alpha^{k}u^{l}$$
where $a_{j} \in Z^{+}$ and $a_{kl}\in Z$ for all $k,l$.
\begin{lem}
  Suppose $F_{j}$ is defined as in (3) for $j\geq 1$. There exists a number $\alpha_{j}>0$ such that $F_{j}(\alpha,u)\geq 0$ for all $\alpha\geq \alpha_{j}$ and $u\geq 0$.
\end{lem}
Proof. Obviously $lim_{\begin{array}{l}
                         \alpha\rightarrow \infty \\
                         u\rightarrow \infty
                       \end{array}}F_{j}(\alpha,u)=\infty$ since $a_{j}>0$. For any $M\geq 0$, there exist $u_{M}>0$ and $\overline{\alpha}_{M}>0$ such that $F_{j}(\alpha,u)\geq M$ whenever $u\geq u_{M}$ and $\alpha \geq \overline{\alpha}_{M}$. Let $M=0$. Then $F_{j}(\alpha,u)\geq 0$ for all $u\geq u_{0}$ and $\alpha\geq \overline{\alpha}_{0}$. For each $u\in [0,u_{0}],\ F_{j}(\alpha,u)$ is a polynomial of $\alpha$. Let it be $p_{u}(\alpha)$. Since $lim_{\alpha\rightarrow \infty}p_{u}(\alpha)=\infty$, there exists $\alpha_{u}>0$ such that $p_{u}(\alpha)\geq 0$ whenever $\alpha\geq \alpha_{u}$. Let $g(u)$ be the infimum of such $\alpha_{u}$. Then $g(u)$ is obviously a continuous function on $[0,u_{0}]$. Let $\beta_{j}=sup_{u\in [0,u_{0}]}g(u)$. Then, for all $\alpha\geq \beta_{j},\ F_{j}(\alpha,u)\geq 0$ if $u\in [0,u_{0}]$. Now, let $\alpha_{j}=max\{ \overline{\alpha}_{0},\beta_{j}\} $. Then $F_{j}(\alpha,u)\geq 0$ for all $\alpha\geq \alpha_{j}$ and $u\geq 0$. \hspace{8cm} \ \ \ \ \ \ $\sharp$
\begin{thm}
  For any positive integer $n$, there is a positive number $M_{n}$ such that for all $\alpha\geq M_{n}$, the function $\Phi$ defined in (1) is strictly positive definite on $R^{d}$ for $1\leq d\leq n$.
\end{thm}
Proof. By Theorem1 of \cite{Dy} and Theorem2.2 of \cite{Mi}, for each $1\leq d\leq n$, if $(-1)^{j}\frac{d^{j}}{dt^{j}}\phi(\sqrt{t})\geq 0$ for $j=0,1,\ldots ,l$ on $(0,\infty)$, where $l=[d/2]+2$, then $\Phi(x)$ is strictly positive definite on $R^{d}$. Let $m_{d}=max\{ \alpha_{1},\ldots ,\alpha_{l}\} $ where $\alpha_{1},\ldots ,\alpha_{l}$ are defined in Lemma2.1. We have $(-1)^{j}\frac{d^{j}}{dt^{j}}\phi(\sqrt{t})\geq 0$ for $j=0,1,\ldots ,l $ on $(0,\infty )$ if $\alpha\geq m_{d}$. Now, let $M_{n}=max\{ m_{1},\ldots , m_{n}\} $. The theorem then follows immediately.
 \hspace{10cm} $\sharp$\\
\\
{\bf Example}. The function defined in (1) is strictly positive definite on $R^{d}$ for $1\leq d\leq 4$ whenever $\alpha \geq 2$.\\
\\
Proof. By Theorem1 of \cite{Dy} and Theorem2.2 of \cite{Mi}, it suffices to show that $(-1)^{j}\frac{d^{j}}{dt^{j}}\phi(\sqrt{t})\geq 0$ for $j=0,1,\ldots ,4$ on $(0,\infty )$, where $\phi$ is as in (2). For $j=0,1,2$, it's trivial. We begin with $j=3$.

For $j=3$,
$$(-1)^{3}\frac{d^{3}}{dt^{3}}\phi(\sqrt{t})=e^{-\alpha(1+u)^{2}}\alpha(\frac{\pi}{2})^{3}(1+u^{2})F_{3}(\alpha,u)$$
where $u=tan\frac{\pi}{2}t$ and
\begin{eqnarray*}
  F_{3}(\alpha,u) & = & 8\alpha^{2}(1+u)^{3}(1+u^{2})^{2}-12\alpha(1+u)(1+u^{2})^{2}\\
                  &   & -24\alpha u(1+u)^{2}(1+u^{2})+8(1+u^{2})u+4(1+2u)(1+u^{2})\\
                  &   & +8(1+u)u^{2}\\
                  & = & 4\alpha(1+u)(1+u^{2})[2\alpha(1+u)^{2}(1+u^{2})-3(1+u^{2})-6u(1+u)]\\
                  &   & +8(1+u^{2})u+4(1+2u)(1+u^{2})+8(1+u)u^{2}.
\end{eqnarray*}
The terms in the brackets equal
\begin{eqnarray*}
  &   & 2\alpha (1+2u+u^{2})(1+u^{2})-3(1+u^{2})-6u(1+u)\\
  & = & 2\alpha (1+u^{2})+2\alpha u^{2}(1+u^{2})-3(1+u^{2})+4\alpha u(1+u^{2})-6u(1+u)\\
  & > & 4\alpha u(1+u^{2})-6u(1+u)+2\alpha u^{2}(1+u^{2})\\
  & = & 4\alpha u+4\alpha u^{3}-6u-6u^{2}+2\alpha u^{2}+2\alpha u^{4}\\
  & = & 2u[2\alpha u^{2}+(\alpha -3)u+2\alpha -3]+2\alpha u^{4}
\end{eqnarray*}
For $3\leq \alpha ,\ (\alpha-3)u+2\alpha-3>0$. For $2\leq \alpha <3$, and $u\geq 1$, it's obvious that $2\alpha u^{2}+(\alpha-3)u+2\alpha-3>0$. For $2\leq \alpha <3$, and $0<u<1$, we have $|(\alpha-3)u|<1$ and $2\alpha-3\geq 1$. This leads to $2\alpha u^{2}+(\alpha-3)u+2\alpha-3>0$. In any case, $F_{3}(\alpha,u)\geq 0$ for $0<u<\infty$ and $\alpha \geq 2$.

We now investigate the case $j=4$. For $j=4$,
\begin{eqnarray*}
  (-1)^{4}\frac{d^{4}}{dt^{4}}\phi(\sqrt{t}) & = & e^{-\alpha (1+u)^{2}}\alpha(\frac{\pi}{2})^{4}(1+u^{2})F_{4}(\alpha,u)\\
                                             & = & e^{-\alpha(1+u)^{2}}\alpha(\frac{\pi}{2})^{4}(1+u^{2})^{2}\left[ \frac{F_{4}(\alpha,u)}{1+u^{2}}\right]
\end{eqnarray*}
where
\begin{eqnarray*}
  \frac{F_{4}(\alpha,u)}{1+u^{2}} & = & 16\alpha^{3}(1+u)^{4}(1+u^{2})^{2}-48\alpha^{2}(1+u)^{2}(1+u^{2})^{2}\\
                                  &   & -96\alpha^{2}(1+u)^{3}u(1+u^{2})  +136\alpha(1+u)u(1+u^{2})\\
                                  &   & +8\alpha(1+u)(1+2u)(1+u^{2})+64\alpha(1+u)^{2}u^{2}\\
                                  &   & +12\alpha(1+u^{2})^{2}+24\alpha(1+u)^{2}(1+3u^{2})-104u^{2}-32u-16\\
                                  &   & -16\cdot \frac{u^{3}(1+u)}{1+u^{2}}.
\end{eqnarray*}
Let $G(\alpha,u):=\frac{F_{4}(\alpha,u)}{1+u^{2}}$. The crucial part of $G(\alpha,u)$ is the first three terms. The three terms equal$$8\alpha^{2}(1+u)^{2}(1+u^{2})\{ (1+u^{2})[2\alpha(1+u)^{2}-6]-12(1+u)u\} . $$
Let $g_{\alpha}(u)=(1+u^{2})[2\alpha(1+u)^{2}-6]-12(1+u)u$. It's easily seen that $g_{\alpha}(2)>0$ for $\alpha\geq 2$. Now,   
 $g_{\alpha}'(u)=2u[2\alpha(1+u)^{2}-6]+(1+u^{2})[4\alpha(1+u)]-12(2u+1)$. For $\alpha=2,\ g_{2}'(u)=8(1+u)(1+u^{2})+8(1+u)^{2}u-36u-12>0$ if $u\geq 2$. Obviously $g_{\alpha}'(u)>0$ for $\alpha>2$ and $u\geq 2$. These lead to $g_{\alpha}(u)>0$ for $\alpha\geq 2$ and $u\geq 2$. It's then easy to see that $G(\alpha,u)>0$ for $\alpha\geq 2$ and $u\geq 2$.

Before analyzing $G(\alpha,u)$ for $\alpha\geq 2$ and $0<u<2$, let's find $\frac{d}{du}G(\alpha,u)$ and $\frac{d}{d\alpha}G(\alpha,u)$ first.
\begin{eqnarray}
  \frac{d}{du}G(\alpha,u) & = & 16\alpha^{3}4(1+u)^{3}(1+u^{2})^{2}+16\alpha^{3}(1+u)^{4}2(1+u^{2})2u\nonumber \\
                          &   & -48\alpha^{2}2(1+u)(1+u^{2})^{2}-48\alpha^{2}(1+u)^{2}2(1+u^{2})2u \nonumber \\                                                                                                                                                                                                                   &   & -96\alpha^{2}3(1+u)^{2}(u+u^{3})-96\alpha^{2}(1+u)^{3}(1+3u^{2})\nonumber \\
                          &   & +136\alpha u(1+u^{2})+136\alpha(1+u)(1+3u^{2})\nonumber \\
                          &   & +8\alpha[8u^{3}+9u^{2}+6u+3]+128\alpha (1+u)u^{2}+64\alpha(1+u)^{2}2u\nonumber \\
                          &   & +24\alpha(1+u^{2})2u+24\alpha 2(1+u)(1+3u^{2})+24\alpha(1+u)^{2}6u\nonumber \\
                          &   & -208u-32-16\cdot \frac{(1+u^{2})(3u^{2}+4u^{3})-u^{3}(1+u)2u}{(1+u^{2})^{2}}\nonumber \\
                          & = & 64\alpha^{3}(1+u)^{3}(1+u^{2})^{2}+64\alpha^{3}(1+u)^{4}(1+u^{2})u\nonumber \\
                          &   & -96\alpha^{2}(1+u)(1+u^{2})^{2}-480\alpha^{2}(1+u)^{2}(u+u^{3})\nonumber \\
                          &   & -96\alpha^{2}(1+u)^{3}(1+3u^{2})+184\alpha u(1+u^{2})\nonumber \\
                          &   & +184\alpha(1+u)(1+3u^{2})+8\alpha[8u^{3}+9u^{2}+6u+3]\nonumber \\
                          &   & +128\alpha(1+u)u^{2}+272\alpha(1+u)^{2}u-208u\nonumber \\
                          &   & -32-16\cdot \frac{(1+u^{2})(3u^{2}+4u^{3})-u^{4}2(1+u)}{(1+u^{2})^{2}}     
\end{eqnarray}
\begin{eqnarray}
  \frac{d}{d\alpha}G(\alpha,u) & = & 48\alpha^{2}(1+u)^{4}(1+u^{2})^{2}-96\alpha(1+u)^{2}(1+u^{2})^{2}\nonumber \\
                               &   & -192\alpha(1+u)^{3}u(1+u^{2})+136(1+u)u(1+u^{2})\nonumber \\                                    &   & +8(1+u)(1+2u)(1+u^{2})+64(1+u)^{2}u^{2}\nonumber \\
                               &   & +12(1+u^{2})^{2}+24(1+u)^{2}(1+3u^{2}) 
\end{eqnarray}
The first two terms of $\frac{d}{d\alpha}G(\alpha,u)$ equal
\begin{eqnarray*}
  &   & 48\alpha^{2}(1+4u+6u^{2}+4u^{3}+u^{4})(1+u^{2})^{2}-96\alpha(1+u)^{2}(1+u^{2})^{2}\\
  & = & 48\alpha^{2}(1+2u+u^{2})(1+u^{2})^{2}+48\alpha^{2}(2u+5u^{2}+4u^{3}+u^{4})(1+u^{2})^{2}\\
  &   & -96\alpha(1+u)^{2}(1+u^{2})^{2}\\
  & \geq &48\alpha^{2}(2u+5u^{2}+4u^{3}+u^{4})(1+u^{2})^{2}\ if\ \alpha \geq 2 \ and\ u\in (0,\infty ).
\end{eqnarray*}
Hence
\begin{eqnarray}
  \frac{d}{d\alpha}G(\alpha,u) & \geq & 48\alpha^{2}(2u+5u^{2}+4u^{3}+u^{4})(1+u^{2})^{2}-192\alpha (1+u)^{3}u(1+u^{2})\nonumber \\
                               &      & +136(1+u)u(1+u^{2})+8(1+u)(1+2u)(1+u^{2})\nonumber \\
                               &      & +64(1+u)^{2}u^{2}+12(1+u^{2})^{2}+24(1+u)^{2}(1+3u^{2})
\end{eqnarray}
The first four terms of (6) equal
\begin{eqnarray}
  &   & (1+u^{2})8[6\alpha^{2}(2u+5u^{2}+4u^{3}+u^{4})(1+u^{2})-24\alpha(1+u)^{3}u\nonumber \\
  &   & +17(1+u)u+(1+u)(1+2u)]\nonumber \\
  & = & (1+u^{2})8[6\alpha^{2}u^{6}+24\alpha^{2}u^{5}+(36\alpha^{2}-24\alpha)u^{4}+(36\alpha^{2}-72\alpha )u^{3}\nonumber \\
  &   & +(30\alpha^{2}-72\alpha+19)u^{2} +(12\alpha^{2}-24\alpha+20)u+1](>0\ if\ u\geq 1)\nonumber \\
  & \geq & (1+u^{2})8[(30\alpha^{2}-72\alpha+19)u^{2}+20u+1]\nonumber \\
  & \geq & (1+u^{2})8[-5u^{2}+20u+1]\nonumber \\
  & \geq & (1+u^{2})8[-5u+20u+1]\ if\ 0<u<1\nonumber \\
  & \geq & (1+u^{2})8[15u+1]\nonumber \\
  & > & 0\ if\ 0<u<1
\end{eqnarray}
Thus
\begin{equation}
 \frac{d}{d\alpha}G(\alpha,u)>0\ if\ \alpha \geq 2\ and\ u\in (0,\infty ).
\end{equation}
Now, let $\alpha=2$. Then
\begin{eqnarray}
  \frac{d}{du}G(2,u) & = & 2^{9}(1+u)^{3}(1+u^{2})^{2}+2^{9}(1+u)^{4}(1+u^{2})u\nonumber \\
                     &   & -3\cdot 2^{7}(1+u)(1+u^{2})^{2} -480\cdot 4(1+u)^{2}u(1+u^{2})\nonumber \\
                     &   & -(2^{8}+2^{7})(1+u)^{3}(1+3u^{2})+u(1+u^{2})[2^{7}+3\cdot2^{6}+2^{5}+2^{4}]\nonumber \\
                     &   & +(1+u)(1+3u^{2})[2^{7}+3\cdot 2^{6}+2^{5}+2^{4}]\nonumber \\
                     &   & +16(3+4u)(1+u^{2})+32(1+3u+2u^{2})u+(1+u)u^{2}\cdot 2^{8}\nonumber \\
                     &   & +(1+u)^{2}u(2^{8}+3\cdot 2^{6}+3\cdot 2^{5})-208u-32\nonumber \\
                     &   & -16\cdot \frac{(1+u^{2})(3u^{2}+4u^{3})-u^{3}(1+u)2u}{(1+u^{2})^{2}}
\end{eqnarray}
We are going to analyze the value of $G(2,u)$ for $1\leq u\leq 2$. Note that for $1\leq u\leq 2,\ (1+u)^{4}\geq 4(1+u)^{2}$. By this the first five terms of (9) can easily be analyzed as follows.
\begin{eqnarray}
  &   & 2^{9}(1+u)^{4}(1+u^{2})u-480\cdot 4\cdot (1+u)^{2}u(1+u^{2})\nonumber \\
  & \geq & 4\cdot 2^{9}(1+u)^{2}(1+u^{2})u-1920(1+u)^{2}u(1+u^{2})\nonumber \\
  & = & 128(1+u)^{2}(1+u^{2})u\ (u\in [1,2])
\end{eqnarray}
Then
\begin{eqnarray}
  &   & 128(1+u)^{2}(1+u^{2})u-3\cdot 2^{7}(1+u)(1+u^{2})^{2}\nonumber \\
  & = & 128[(1+u^{2})^{2}u+2(1+u^{2})u^{2}]-3\cdot 2^{7}(1+u)(1+u^{2})^{2}\nonumber \\
  & = & -256(1+u^{2})^{2}u-384(1+u^{2})^{2}+256(1+u^{2})u^{2}.
\end{eqnarray}
Now,
\begin{eqnarray}
  &   & 2^{9}(1+u)^{3}(1+u^{2})^{2}-(2^{8}+2^{7})(1+u)^{3}(1+3u^{2})\nonumber \\
  & = & 2^{9}(1+u)^{3}(1+u^{2})^{2}-(2^{8}+2^{8})(1+u)^{3}(1+3u^{2})+2^{7}(1+u)^{3}(1+3u^{2})\nonumber \\
  & = & 2^{9}(1+u)^{3}[1+2u^{2}+u^{4}-1-3u^{2}]+2^{7}(1+u)^{3}(1+3u^{2})\nonumber \\
  & = & 2^{9}(1+u)^{3}[u^{4}-u^{2}]+2^{7}(1+u)^{3}(1+3u^{2})\nonumber \\
  & \geq & 2^{7}(1+u)^{3}(1+3u^{2}).\ (u\in [1,2])
\end{eqnarray}
The sum of (11) and (12) is $128u^{5}+1024u^{4}+256u^{2}+128u-256$ which is larger than $0$ if $u\in [1,2]$. Thus, the sum of the first five terms of (9) is positive. It's trivial to see that the sum of the other terms is also positive. Hence
\begin{equation}
  \frac{d}{du}G(2,u)>0
\end{equation}
for $u\in [1,2]$. Since $G(2,1)>0$, by (8) and (13), we get $G(\alpha,u)\geq 0$ for $\alpha \geq 2$ and $u\in [1,2]$.

The next task is to analyze $G(\alpha,u)$ for $\alpha \geq 2$ and $u\in (0,1)$. We first let $\alpha =2$ and $u\in (0,1)$. In this case,
\begin{eqnarray*}
  G(2,u) & = & 16\cdot 2^{3}(1+u)^{4}(1+u^{2})^{2}-48\cdot 2^{2}(1+u)^{2}(1+u^{2})^{2}\\
         &   & -96\cdot 4 (1+u)^{3}(u+u^{3})+136\cdot 2(u+u^{2})(1+u^{2})\\
         &   & +16(1+u)(1+2u)(1+u^{2})+128(1+u)^{2}u^{2}\\
         &   & +24(1+u^{2})^{2}+48(1+u)^{2}(1+3u^{2})-104u^{2}\\ 
         &   & -32u-16-16\cdot \frac{u^{3}(1+u)}{1+u^{2}}\\
         & = & 2^{7}u^{8}+2^{9}u^{7}+(5\cdot 2^{6}+2^{7})u^{6}+280u^{4}+96u^{3}-120u^{2}\\
         &   & +128u+8-16\cdot \frac{u^{3}(1+u)}{1+u^{2}}\\
         & \geq & 8u+8-16(u^{3}+u^{4})+96u^{3}+280u^{4}\ (since\ u\in [0,1])\\
         & = & 8u+8+80u^{3}+264u^{4}\\
         & \geq & 8. 
\end{eqnarray*}
This result together with (8) shows that $G(\alpha,u)\geq 0$ for $\alpha\geq 2$ and $u\in (0,1)$.

Our conclusion thus follows. \hspace{7cm} \ \ \ $\sharp$\\
\\
{\bf Remark.} For dimensions higher than four, the value of $\alpha$ can be analyzed analogously. It's a bit more tedious, but not difficult. For each increase of the dimension by two, the order of the derivative of $\phi(\sqrt{t})$ has to be increased by one.\\

\parbox[t]{2.5in}{Lin-Tian Luh\\ Institut fuer Numerische und\\ Angewandte Mathematik\\ Universitaet Goettingen\\ Lotze Str. 16-18\\ 37083 Goettingen\\ Germany\\ luh@math.uni-goettingen.de\\ (in July and August)}\hspace{1.5cm}
\parbox[t]{3.5in}{Lin-Tian Luh\\ Department of Mathematics\\ Providence University\\ Shalu Town\\ Taichung County\\ Taiwan\\ ltluh@pu.edu.tw\\ (in normal time except July\\ and August)}


\begin{thebibliography}{99}
\bibitem{Dy}N. Dyn,
{\em Interpolation and Approximation by Radial and Related Functions,}
Approximation Theory VI,(C.K. Chui, L.L. Schumaker and J. Ward eds.), Academic Press,(1989),211-234.

\bibitem{Mi}C.A. Micchelli,
{\em Interpolation of Scattered Data:Distance Matrices and Conditionally Positive Definite Functions,}
Constr. Approx. 2(1986),11-22.

\end{thebibliography}
\end{document}